%% file: knight.tex
\documentclass[12pt]{article}
\usepackage{latexsym}
\usepackage{fullpage,epsfig}
\usepackage{fullpage}
\usepackage{amssymb}

\newcommand{\epsf}[1]{\epsfbox{#1}}

\def\qed{$\hfill{\vrule height 3pt width 5pt depth 2pt}$}

\newfont{\bbold}{msbm10 scaled \magstep1}
\newfont{\bbolds}{msbm7 scaled \magstep1}

\newcommand{\ns}{\mbox{\bbold N}}
\newcommand{\zs}{\mbox{\bbold Z}}
\newcommand{\zss}{\mbox{\bbolds Z}}
\newcommand{\qs}{\mbox{\bbold Q}}
\newcommand{\rs}{\mbox{\bbold R}}
\newcommand{\cs}{\mbox{\bbold C}}

\newcommand{\bm}[1]{\mbox{\boldmath \ensuremath{#1}}}
\newcommand{\bs}[1]{\mbox{\boldmath \ensuremath{\scriptstyle #1}}}

\newcommand{\bx}{\bar x}
\newcommand{\by}{\bar y}
\newcommand{\bu}{\bar u}

\newcommand{\GL}{\mathbb{L}}

\newcommand{\Ref}[1]{(\ref{#1})}
\newcommand{\beq}{\begin{equation}}
\newcommand{\eeq}{\end{equation}}
\newcommand{\gf}{generating function}
\newcommand{\gfs}{generating functions}

\def\cqfd{\par\nopagebreak\rightline{\vrule height 3pt width 5pt depth 2pt}
\medbreak}

 \newtheorem{Theorem}{Theorem}
 \newtheorem{Proposition}[Theorem]{Proposition}

\newtheorem{Lemma}[Theorem]{Lemma}

\title{Walks confined in a quadrant are not always D-finite}
\author{
\parbox{8cm}{
 {\sc Mireille Bousquet-M{\'e}lou\thanks{Partially
supported by the INRIA, via the cooperative research action
\textsf{Alcophys}.}}\\
{\small CNRS, LaBRI, Universit\'e Bordeaux 1 \\
351 cours de la Lib\'eration \\
33405 Talence Cedex,  France\\
{\tt mireille.bousquet@labri.fr }
}}
\parbox{7cm}{
{\sc Marko Petkov\v sek\thanks{Partially supported by 
MZT RS under grant P0-0511-0101.}} \\
{\small Department of Mathematics \\
University of Ljubljana \\
Jadranska 19, SI-1000 Ljubljana, Slovenia \\
 {\tt Marko.Petkovsek@fmf.uni-lj.si}}
}}

\date{}

\begin{document}
\maketitle

\begin{abstract}
We consider planar lattice walks that start from a prescribed
position, take their steps  in a given finite subset of $\zs ^2$, and always
stay in the quadrant $x\ge 0, y\ge 0$. We first give a criterion
which guarantees that the length \gf    \ of these walks is D-finite,
that is, satisfies a linear differential equation with polynomial
coefficients. This criterion applies, among others, to the ordinary
square lattice walks. Then, we prove that walks that start from
$(1,1)$, take their steps in $\{(2,-1), (-1,2)\}$ and stay in the
first quadrant have a non-D-finite \gf . Our proof relies on a
functional equation satisfied by this \gf , and on elementary complex
analysis.
\end{abstract}

\section{Introduction}
The enumeration of lattice walks is one of the most venerable topics
in enumerative combinatorics, which has numerous applications in
probabilities~\cite{feller,lawler,spitzer}. These walks take their steps in a
finite subset $\frak S$ of $\zs ^d$, and might be constrained in
various ways. One can only cite a small percentage of the relevant
litterature, which dates back at least to the next-to-last century
~\cite{andre,motzkin,gessel-proba,kreweras,mohanty,narayana}.
Many recent publications show that the topic is still
active~\cite{firenze,mbm-slitplane,prep-GS,gessel-zeilberger,grabiner,niederhausen1,niederhausen2}.

After the solution of many explicit problems, certain patterns have
emerged, and a more recent trend consists in 
developing methods that
are valid for generic sets of steps.
Special attention is being paid
to the {\em nature\/} of the \gf \ of the walks under
consideration. For instance, the \gf \ for unconstrained walks on the
line $\zs$ is rational, while the \gf \ for walks constrained to stay
in the half-line $\ns$ is always algebraic~\cite{banderier-flajolet}.
This result has often been described in terms of {\em partially directed\/}
2-dimensional walks confined in a quadrant (or {\em generalized Dyck
walks\/}~\cite{duchon,gessel-factor,lab3,lab2,merlini}), 
but is, essentially, of a 1-dimensional nature.

Similar questions can be addressed for  {\em real\/} 2-dimensional
walks. Again, the \gf \ for unconstrained walks starting 
from a given 
point is clearly rational. Moreover, the argument used for
1-dimensional walks confined in $\ns$ can
be 
recycled to prove that the \gf \ for the  walks that
stay in the half-plane $x\ge 0$ is always algebraic. What about
doubly-restricted walks, that is, walks that are confined in the
quadrant $x\ge 0, y \ge 0$? It would be satisfactory if the hierarchy
unconstrained walks / mono-constrained walks / bi-constrained walks
could match the classical hierarchy of \gfs : rational series,
algebraic series, D-finite series\footnote{A series $F(t)$ is
D-finite if it satisfies a linear differential equation with 
polynomial coefficients in $t$.} (also called holonomic series). 
A rapid inspection of the most standard 
cases only corroborates this hope. For instance, the \gf \ for walks on
the square lattice (with North, East, South and West steps) 
that start from the origin and stay in the first
quadrant is
$$\sum_{m, n \ge 0} {{m+n} \choose m} {m \choose {\lfloor m/2\rfloor}}
{n \choose {\lfloor n/2\rfloor}} t^{m+n}
=
\sum_{ n \ge 0} {n \choose {\lfloor n/2\rfloor}} 
{{n+1} \choose {\lceil n/2\rceil}} t^{n},$$
which is a D-finite series.  The first expression comes from the fact that
these walks are  shuffles of two prefixes of Dyck walks, and the
Chu-Vandermonde identity transforms it into the second simpler expression.

\begin{figure}[ht]
\begin{center}
\epsf{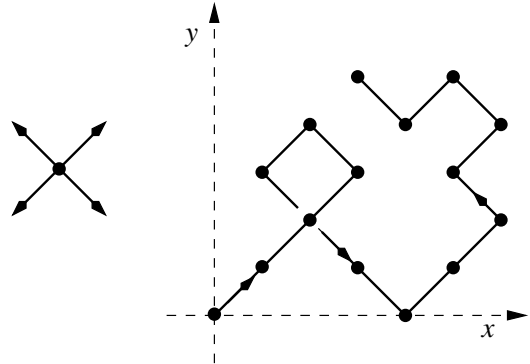}
\end{center}

\caption{A walk on the diagonal square lattice confined in the first quadrant.}
\label{diagonal}
\end{figure}

The case of the diagonal square lattice, where the steps are
North-East, South-East, North-West and South-West (Figure~\ref{diagonal}) is
even simpler: by projecting the walks on the $x$- and $y$-axes, we
obtain two decoupled prefixes of Dyck paths, so that the \gf \ for
walks  in the first quadrant is now
$$\sum_{n \ge 0} {n \choose {\lfloor n/2\rfloor}}^2 t^n,$$
another D-finite series.
In both cases, the number of $n$-step walks can be shown to grow
asymptotically like $4^n/n$, which prevents the corresponding \gf \
from being algebraic (see~\cite{flajolet} for the possible asymptotic
behaviours of coefficients of algebraic series).

In Section \ref{section-sufficient}  of this paper, we shall
generalize this result by 
proving that, if the set of steps $\frak S$ is symmetric with respect to the
$x$-axis and satisfies a {\em small height variation\/} condition,
then the \gf \ for walks with steps in $\frak S$, starting from any
given point $(i_0,j_0)$, is D-finite. This result covers the above two
cases.

However, and most importantly, we shall also prove in
Section~\ref{section-knight} that this holonomy result does not hold
for {\em any\/} set of 
steps: walks that start from $(1,1)$, take their steps in ${\frak S}=
\{(2,-1), (-1,2)\}$ and always stay in the first quadrant have a
non-D-finite \gf . The central point of our proof is the study of a
series $G(x)$, defined by an equation of the form
$$G(x)+G(\xi(x))=A(x),$$
where $\xi(x)$ and $A(x)$ are explicit algebraic series in $x$. We
consider  the solution $G(x)$ to this equation as a function of a
complex variable $x$, and prove that it has {\em
 infinitely many singularities\/}, 
which prevents it from
being D-finite. Hence our proof is based on complex analysis.
To our knowledge, there is no
classification of the solutions to this 
type of equation.
In some very
specific cases (like $\xi(x)=x^p$ or $\xi(x)=cx$) some
hypertranscendence results\footnote{A series $F(x)$  
is hypertranscendental if it does not satisfy any polynomial
differential equation of the type $P(x,F(x), F'(x), \ldots ,
F^{(k)}(x))=0$, where $P$ is a polynomial.} 
have been obtained, either by some {\em ad hoc\/}
methods~\cite{ishizaki,loxton,nishioka}, 
or  via general results about ``very'' lacunary
series~\cite{lipshitz-rubel}.

These two sections raise the question of a classification of the sets
$\frak S$ according to the nature of the \gf \ for walks in a quadrant
that take their steps in $\frak S$. Let us mention that some sets of
steps, like $\frak S= \{(1,1), (0,-1), (-1,0)\}$ yield, for non-trivial
reasons, algebraic \gfs \
\cite{mbm-kreweras-conf,mbm-kreweras,gessel-proba,kreweras,niederhausen1}.

Finally, in Section~\ref{section-recurrences}, we say a few words
about the closely related topic of multidimensional linear recurrences
with constant coefficients, and prove a non-holonomy result that was
announced (but not proven) in~\cite{bousquet-petkovsek}.

\bigskip
Let us conclude this introduction with a few more formal definitions
on walks and power series.

Let $\frak S$ be a finite subset of $\zs ^2$. A walk with steps in
$\frak S$ is a
finite sequence $w = (w_0, w_1, \ldots , w_n)$  of vertices of $\zs
^2$ such that $w_{i}-w_{i-1} \in \frak S$ for $1 \le i \le
n$. The number of 
steps, $n$, is the {\em length\/} of $w$. The  starting point
of $w$ is $w_0$, and its  endpoint is
$w_n$.  The {\em complete \gf\/} for  a set $\cal A$ of walks starting
from a given point $w_0$
is the series
$$ A(x,y;t)=\sum_{n \ge 0}t^n \sum_{i,j \in \zss} a_{i,j}(n) x^i 
y^j ,$$ where $a_{i,j}(n)$ is the number  of walks of  $\cal A$
that have length $n$  and  end at $(i,j)$. This series is a formal
power series in $t$ whose coefficients are polynomials in $x, y, 1/x,
1/y$. We shall often denote $\bx =1/x$ and $\by =1/y$.
The length \gf \ for walks of $A$ is simply
$$A(t)= \sum_{n \ge 0} a(n) t^n$$
 where $a(n)$ is the number  of walks of  $\cal A$
that have length $n$.
Note that $A(t)=A(1,1;t)$.

 Given a ring $\GL$ and $k$ indeterminates
$x_1, \ldots , x_k$, we denote by
 $\GL[x_1, \ldots , x_k]$ (resp.~$\GL[[x_1, \ldots , x_k]]$) the ring of polynomials 
(resp.~formal power series)
in $x_1, \ldots , x_k$ with coefficients in $\GL$.
%
%
%
 If $\GL$ is a field, we denote by
 $\GL(x_1, \ldots , x_k)$ the field of rational functions in $x_1,
\ldots , x_k$ with coefficients in $\GL$.
%

Assume $\GL$ is a field. A series $F$ in $\GL[[x_1, \ldots , x_k]]$ is
{\em rational\/} if 
there exist  polynomials $P$ and $Q$
in $\GL[x_1, \ldots ,
x_k]$, with $Q\not = 0$, such that $QF=P$. 
It is {\em algebraic\/}
(over the field  $\GL(x_1, \ldots , x_k)$) if
there exists a non-trivial polynomial $P$ with coefficients in
$\GL$ such that 
$P(F,x_1, \ldots , x_k)=0.$ The sum and product of algebraic series
is algebraic.  
%

The series $F$  is {\em D-finite\/} (or {\em holonomic\/})  
if the partial derivatives of $F$ span a finite 
dimensional vector space over the field $\GL(x_1, \ldots , x_k)$
(this vector space is a subspace  of the fraction field of $\GL[[x_1, \ldots , x_k]]$);
see~\cite{stanleyDF} for the one-variable case,
and~\cite{lipshitz-diag,lipshitz-df} otherwise.  In other 
words, for $1\le 
i\le k$, the series $F$ satisfies a non-trivial partial differential
equation of the form
$$\sum_{\ell=0}^{d_i}P_{\ell,i}\frac{\partial ^\ell
F}{\partial x_i^\ell} =0,$$
where $P_{\ell,i}$ is a polynomial in the $x_j$.
Any algebraic series is holonomic. The
sum and product of two 
holonomic series is still holonomic. The
specializations of a holonomic series (obtained by giving 
values from $\GL$  to some of the variables) are holonomic, if
well-defined. Moreover, if $F$ is an {\em algebraic\/} series and 
$G(t)$ is a holonomic series of one variable, then the substitution
$G(F)$ (if well-defined) is
holonomic~\cite[Prop.~2.3]{lipshitz-df}. 

\section{A sufficient condition for holonomy}\label{section-sufficient}
Let $\frak S$ be a finite subset of $\zs^2$. We say that  $\frak S$ is
symmetric with respect to the $x$-axis if 
$$(i,j) \in {\frak S} \Rightarrow (i,-j) \in {\frak S}.$$
We say that $\frak S$ has small height variation if 
$$(i,j) \in {\frak S} \Rightarrow |j|\le 1.$$
The usual square lattice steps satisfy these two conditions. So do the
steps of the diagonal square lattice (Figure~\ref{diagonal}).

\begin{Theorem}\label{theorem-sufficient}
Let ${\frak S}$   be a finite subset of $\zs^2$ that is symmetric with
respect to the $x$-axis and has small height variations. Let
$(i_0,j_0) \in \ns^2$. Then the length \gf \ for walks that start from
$(i_0, j_0)$, take their steps in $\frak S$ and stay  in the first
quadrant is D-finite.
\end{Theorem}
We shall need the following preliminary result, which does not require
any property on $\frak S$.

\begin{Proposition}\label{propo:alg}
Let ${\frak S}$   be a finite subset of $\zs^2$.  Let
$(i_0,j_0) \in \ns \times \zs $. Then the complete \gf \ for walks
that start from $(i_0, j_0)$, take their steps in $\frak S$ and stay in
the right half-plane $x\ge 0$ is algebraic. 
\end{Proposition}
{\bf Proof.} 
This result  is basically of a 1-dimensional nature: projecting the
walks on the $x$-axis reduces it to the 
enumeration of walks on the half-line $\ns$, starting
from $i_0$, in which each step of size $i$ is weighted 
by a Laurent polynomial in $y$:
$$\sum_{j : (i,j) \in \frak S} y^j.$$
The weight of a walk is taken to be the product of the weights of its
steps. The material of~\cite{banderier-flajolet,bousquet-petkovsek}
or~\cite{gessel-factor} can be readily extended to weighted 
paths, where the weights belong to any field of characteristic $0$
(here, $\qs(y)$), and we conclude that the
complete \gf \ for 
walks in the right-half plane is algebraic over $\qs(x,y,t)$.
\cqfd


\noindent{\bf Note.} The above proposition  is
 very close  to the results on partially directed 
2-dimensional walks mentioned in the introduction. However, in most
references,  the authors only consider the case where all steps  are
\emph{directed}; that is,  $i \ge 0$ for all step $(i,j) \in 
{\frak S}$.

\bigskip

\noindent
{\bf Proof of Theorem~\ref{theorem-sufficient}.} 
Let $\cal Q$ denote the set of  walks that start from
$(i_0, j_0)$, take their steps in $\frak S$ and stay  in the first
quadrant. 
We shall prove that walks of $\cal Q$ are, roughly speaking,
``equivalent to'' walks
in the right half-plane  ending on the $x$-axis.

We claim that it suffices to prove that the subset of $\cal
Q$ consisting of
the walks that hit the $x$-axis at some point has a
D-finite \gf . Indeed, the remaining walks are, by a vertical
translation, in one-to-one correspondence with walks that start from
$(i_0,j_0-k)$, for some $k \in [1,j_0]$, stay in the first quadrant
and hit the $x$-axis.

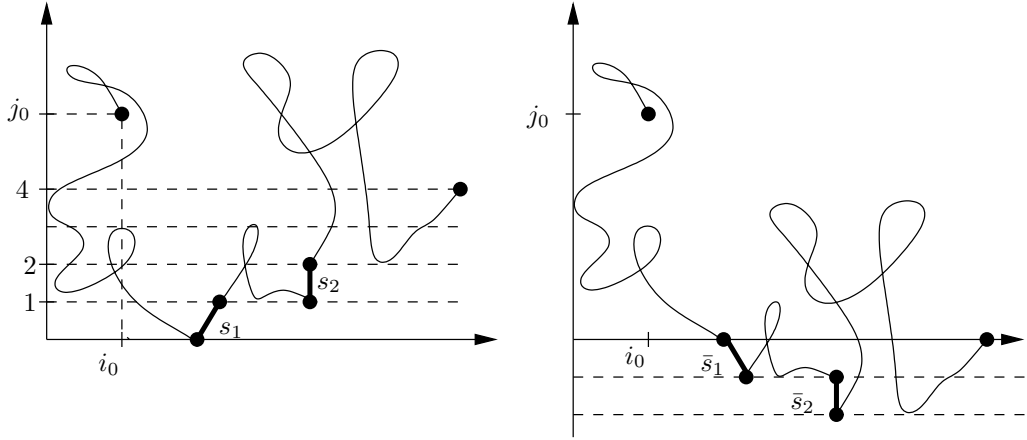
\begin{figure}[ht]
\begin{center}
\input{bijection.pstex_t}
\end{center}
\caption{The bijection between some walks in the quadrant and walks in the
right half-plane ending on the $x$-axis.}
\label{bijection}
\end{figure}

Let us first focus on the set ${\cal Q}^{(0)}_e$ of walks in the first
quadrant that start from
$(i_0,j_0)$,  hit the $x$-axis and end at an even ordinate. These
walks are in bijection with the set ${\cal H}_0$ of walks that start
from $(i_0,j_0)$, stay in the right half-plane and end on the
$x$-axis. This bijection, illustrated by Figure~\ref{bijection}, is a
mere adaptation of a classical 1-dimensional
correspondence, 
which establishes that Dyck prefixes ending at an even
ordinate are equivalent to bilateral Dyck walks (see the Catalan
factorisation in~\cite{cori}.) 
Starting from a walk $w$ of ${\cal
Q}^{(0)}_e$, ending at level $2k$, we denote by $s_1, \ldots , s_k$
the steps that follow the last visit of $w$ to a point of level
(ordinate) $0,\ldots, k-1$. Replacing these $k$ steps by their
symmetric steps
with respect to the $x$-axis yields a walk $\bar w$ that
belongs to  ${\cal H}_0$. The ordinate of the lowest point(s) visited
by $\bar w$ is $-k$. 
Conversely, the steps  $\bar s_1, \ldots , \bar
s_k$ of $\bar w$ that we have to
flip back to recover $w$ are 
the first steps of $w$ that lead
to level $-j$, for $1 \le j \le k$. 

Let $H(x,y;t)$ be the complete \gf \ for  walks in the
right half-plane, starting from $(i_0,j_0)$.  
By Proposition~\ref{propo:alg},
this series is algebraic, hence D-finite. 
The length \gf \ for
walks of ${\cal H}_0$ is 
obtained by extracting the 
coefficient of $y^0$ in the \gf \ $H(1,y;t)$. 
Extracting the constant term  of a D-finite series is known to give
another D-finite series~\cite{lipshitz-diag}: thus the \gf \ for
walks of ${\cal H}_0$, and hence of 
${\cal Q}_e^{(0)}$, is D-finite.

A similar argument holds for the set  ${\cal Q}^{(0)}_o$ of walks in the first
quadrant that start from
$(i_0,j_0)$,  hit the $x$-axis and end at an odd ordinate: they are in
one-to-one correspondence with  the set ${\cal H}_{-1}$ of walks that start
from $(i_0,j_0)$, stay in the right half-plane and end at level $-1$.
 The \gf \ for walks of  ${\cal H}_{-1}$ is the coefficient
of $y^0$ in $yH(1,y;t)$,
 and hence is D-finite.
Given that the sum of D-finite
series is D-finite, this concludes the proof of
Theorem~\ref{theorem-sufficient}. \cqfd

\noindent
{\bf Remark.} A slightly stronger result is proved
in~\cite{mbm-kreweras-conf} via a 
functional equation approach: under the assumptions of
Theorem~\ref{theorem-sufficient},  the {\em complete\/} \gf \ $Q(x,y;t)$
for walks in  the quadrant, counted by their length and coordinates of the
endpoint, is also D-finite. 

\section{The knight walk is not holonomic}\label{section-knight}

\subsection{The main result}
\begin{figure}[ht]
\begin{center}
\epsf{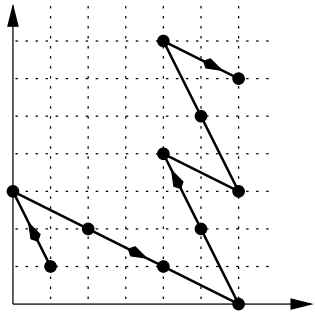}
\end{center}

\caption{A knight walk.}
\label{knight}
\end{figure}

 We study walks  that start at $(1,1)$, take their steps in
 $\{ (-1,2), (2,-1)\}$ and  stay in the first quadrant. We call them
 {\em knight walks\/}, since their steps correspond to two of the knight
 moves on a chessboard (Figure~\ref{knight}). We note that a walk
 ending at $(i,j)$ always has 
  $i+j-2$ steps: hence the information contained in the complete
 \gf \ is actually already contained in the following bivariate series:
 $$Q(x,y)= \sum_{i \ge 0, j\ge 0} Q_{i,j} \, x^i y^j,$$
where $Q_{i,j}$ denotes the number of knight walks ending at  $(i,j)$.

The coefficients $ Q_{i,j}$ satisfy the following recurrence relation:
\beq Q_{i,j}= \left\{ 
\begin{array}{ll}
    0  & \hbox{ if } i<0  \hbox{ or }   j<0 ,\\
    1   & \hbox{ if } i=j=1, \\
    Q_{i+1,j-2}+Q_{i-2,j+1} & \hbox{ otherwise. }
\end{array} \right. \label{rec1} \eeq
This recurrence allows us to compute the numbers  $ Q_{i,j}$
inductively, for instance diagonal by diagonal. The first few values
are given in Table~\ref{table-knight} below, in which the zero entries are left
out. The non-zero entries lie on the
lines  $i=j \hbox{ mod } 3$.

\begin{table}[ht]
\begin{small}
$$\begin{array} {c|ccccccccccccccc}
j &&&&&&& \\
\uparrow &&&&&&& \\

12 &24 \\

11&  &  & 108 & \\

10&   & 24 &  &  & 312  &  &  &  &  &  &  \\

9& 6  &  &  & 84 &  &  & 720  &  \\

8&  &  & 24  &  &  & 204   & & & 1440 \\

7 &     & 6 &  &  & 60  &  &  & 408   &  \\

6 & 2  &  &  & 18 &  &  & 120 &  &  & 720  &  \\

5 &  &  & 6  &  &  &36  &  &  &204  &  \\

4 &  & 2 &  &  &12  &  &  &60  &  &  &312  \\

3 & 1  &  &  & 4 &  &  & 18 &  &  &84  \\

2&  &  & 2 &  &  & 6 &  &  &24  &  &  & 108 \\

1& & 1  &  &  & 2 &  &  &6  &  &  & 24  \\
0&  &  &  & 1  &  &  & 2 &  &  & 6 &  &  & 24  \\ \hline 

& 0 & 1 & 2 & 3 & 4 & 5 & 6 & 7 & 8 & 9 & 10 & 11 & 12 & \rightarrow \ i
\end{array}$$
\end{small}

\caption{The number $Q_{i,j}$ of knight walks ending at $(i,j)$.}
\label{table-knight}\end{table}

\begin{figure}[ht]
\begin{center}
\epsfxsize=6cm \epsfysize=6cm
\epsf{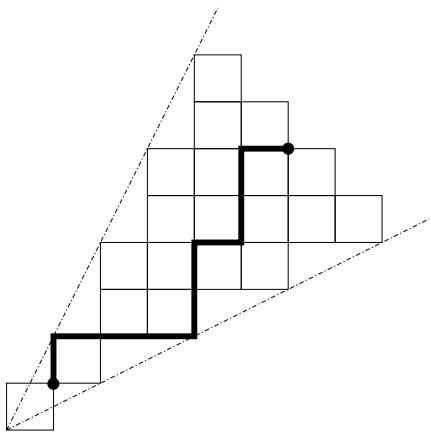}
\end{center}
\caption{Another description of the knight walk.}
\label{redresse}
\end{figure}

Applying the transformation $(i,j)\rightarrow ((2i+j)/3,(i+2j)/3)$
shows that $Q_{i,j}$ is also the number of walks made of North and
East steps, that start from $(1,1)$, end at $((2i+j)/3,(i+2j)/3)$
and always stay above the line $2y=x$ and below the line $2x=y$ (see
Figure~\ref{redresse}). 
One thus recognizes a problem already studied
by  Gessel~\cite[Section~4]{gessel-proba}. But the results presented here are
stronger than those of~\cite{gessel-proba}.  
Ignoring the two
boundary lines  gives the  following simple bound: 
$$Q_{i,j}\le {{i+j-2}\choose {(2i+j-3)/3}}.$$
In particular,
\beq Q_{3i,0} \le {{3i-2}\choose {i-1}}.\label{bound}\eeq
This bound could be sharpened by counting walks that stay above the
line $2y=x$ (see e.g.~\cite{motzkin} or~\cite[Example~4]{bousquet-petkovsek}),
but the above bound will be enough for our purpose.

By summing the
recurrence relation we obtain:
\begin{eqnarray*}
Q(x,y)=\sum_{i,j \ge 0} Q_{i,j}x^i y^j &=& xy
+ \sum_{i,j \ge 0} (Q_{i+1,j-2} +
Q_{i-2,j+1}) x^i y^j\\
&=& xy
+ y^2/x (Q(x,y)-Q(0,y)) + x^2/y (Q(x,y)-Q(x,0)) ,
\end{eqnarray*}
that is
\beq ( xy-x^3-y^3)Q(x,y)=x^2y^2-G(x)-G(y)\label{equation-knight} \eeq
where 
\beq    \label{G-def}
G(x)=x^3 \sum_{i\ge 0} Q_{i,0} x^{i}=x^3Q(x,0)
\eeq
counts knight walks ending on the $x$-axis.
(We have used the symmetry of the problem in $x$ and $y$.)
Note that the length \gf \ for all knight walks is $t^{-2}Q(t,t)$. The
above equation, combined with the
elementary properties of D-finite series, imply
that the following
three statements are equivalent:

--- the bivariate \gf \ $Q(x,y)$ for knight walks is D-finite,

--- the \gf \ $Q(x,0)$ for knight walks ending on the $x$-axis is
    D-finite,

--- the length \gf \ $t^{-2}Q(t,t)$ for knight walks is    D-finite.

\noindent Our main result asserts that none of these statements hold.
\begin{Theorem}\label{theorem-knight}
The length \gf \ for walks that start from $(1,1)$, take their steps
in  $\{ (-1,2), (2,-1)\}$ and always stay in the first quadrant 
(knight walks) is
not D-finite. Nor is the \gf \ for knight walks that end on the $x$-axis.
\end{Theorem}
Looking at Table~\ref{table-knight}, one might still have some hope
that the numbers $Q_{i,j}$ are not so bad. In particular, they seem to
have small prime factors. This pattern actually does not go on, and,
in case there would still be a doubt, the following proposition  reinforces the
non-holonomic character of these numbers.
\begin{Proposition}\label{propo-diagonal}
The length \gf \ for knight walks ending on the main diagonal $x=y$ is
not D-finite.
\end{Proposition}

\subsection{The kernel method}
The so-called {\em kernel method\/}  solves completely the
functional equation~\Ref{equation-knight}. This method has been around
since, at least, the 
70's, and is currently the subject of a certain rebirth (see the
references
in~\cite{hexacephale,banderier-flajolet,bousquet-petkovsek}).

Applied to our equation, this method consists in coupling the
variables $x$ and $y$ so as to cancel the {\em kernel\/} $xy-x^3-y^3$;
this yields the 
missing information about
the series $G(x)$.
More precisely, let $\xi(x)$ be the unique formal power series in $x$
satisfying $$x\xi-x^3-\xi^3=0.$$
The Lagrange inversion formula, applied to  $\xi(x)/x$, provides 
an explicit expression for $\xi(x)$:
\beq \xi(x)=x^2 \sum_{m \ge 0} \frac{x^{3m}}{2m+1} {3m \choose m}=
O(x^2).\label{xi-def}\eeq
Replacing $y$ by $\xi(x)$ in \Ref{equation-knight} gives a functional
equation that defines the power series $G(x)$:
\beq G(x)+G(\xi(x))=x^2\xi(x)^2.\label{main} \eeq
Indeed we obtain, after iterating this equation infinitely
many times:
\beq \label{eq:iterated}
G(x)=\sum_{i\ge 0} (-1)^i  \left( \xi^{(i)}(x) \xi^{(i+1)}(x)
\right) ^2
\eeq
where $\xi^{(i)} = \xi \circ \cdots \circ \xi $ is the $i$th
iterate of $\xi$. Note that $\xi^{(i)}(x)=O(x^{2^i})$, so that
%
%
the sum is convergent in the ring $\cs [[x]]$. 
Replacing $G(x)$ by the above explicit value in
\Ref{equation-knight} would give an expression for $Q(x,y)$.

However, we shall not exploit these expressions, but rather the functional
equation~\Ref{main}, to prove that  $G(x)$, hence $Q(x,y)$, is not
D-finite. 
Our  proof will be  of an  analytic nature.
The idea is to consider $G(x)$ as a  function of a complex variable
$x$ and study its singularities. Using functional equations like
\Ref{main}, we shall
 build new singularities of $G$ from old ones -- and end up with infinitely
many singularities, thus proving that $G$ cannot be holonomic.

However, even though \Ref{main}  defines $G(x)$  uniquely,
this  equation itself is not sufficient for our 
 purpose: we shall introduce the other  two
 roots $\xi_1(x)$ and $\xi_2(x)$ of the kernel, and the
corresponding analogues of~\Ref{main}, to obtain enough singularities.

\subsection{The roots of $x^3+y^3=xy$}

As a polynomial in $y$, the kernel $xy-x^3-y^3$ has three roots. Only
one of them, given by~\Ref{xi-def}, is a power series in $x$. We shall
denote it, from now on, by $\xi_0$:
\begin{eqnarray*}
\xi_0(x)&=& x^2+x^5+3x^8+12x^{11}+55x^{14}+273x^{17}+ \cdots
\end{eqnarray*}
The other two roots are power series in $\sqrt{x}$, and their
expansion can be computed inductively:
\begin{eqnarray*}
\xi_1(x)&=&\displaystyle +\sqrt{x}-\frac{x^2}{2} - \frac{3}{8} x^3\sqrt{x}
 -\frac{x^5}{2} -\frac{105}{128} x^6\sqrt {x} - \cdots \\
 & & \\
\xi_2(x)& =& \displaystyle -\sqrt{x}-\frac{x^2}{2} +\frac{3}{8} x^3\sqrt{x}
 -\frac{x^5}{2} +\frac{105}{128} x^6\sqrt {x} - \cdots 
\end{eqnarray*}
Of course, $\xi_2(x)$ is derived from $\xi_1(x)$ by replacing $\sqrt
x$ by $-\sqrt x$. Guided by the above expressions, let us write
\begin{eqnarray}
\xi_1(x) & =& +\sqrt x \, \psi(x)-\phi(x) \label{xi1-phi-psi}\\
\xi_2(x) & =& -\sqrt x \, \psi(x)-\phi(x) \nonumber
\end{eqnarray}
where $\psi$ and $\phi$ are formal power series in $x$. As the  three
roots sum to zero, one has 
$$ \phi(x)=\xi_0(x)/2. $$
In order to compute the coefficients of $\psi$, we shall use again the
Lagrange inversion formula (LIF). 
Let $\chi(x)$ be defined by
$$ \xi_1(x)=\frac{\sqrt x}{\sqrt{1+\chi(x)}}. $$
Then
$$\chi={x\sqrt x} \left({1+\chi}\right)^{3/2},$$
so that the LIF gives
$$
\xi_1(x)=\sqrt x \left( 1 - \sum_{n\ge 1} \frac{(x\sqrt x)^n }{2n}
{{3(n-1)/2} \choose {n-1}}\right).
$$
Separating the even/odd values of $n$ gives explicit expressions of
$\psi(x)$ and $\phi(x)$. 
The following lemma summarizes the results thus obtained.
\begin{Lemma}\label{three-roots}
Let $\xi$ and $\psi$ be the following power series in $x$:
\begin{eqnarray*}
 \xi&=& x^2 \sum_{m \ge 0} \frac{x^{3m}}{2m+1} {3m \choose m},\\
\psi&=&1 -\sum_{m \ge 1} 
\frac{m! (6m)! }{(6m-1)(2m)!^2(3m)!} \frac{x^{3m}}{16^m}.
\end{eqnarray*}
Then the three roots of $x^3+y^3=xy$ are 
\begin{eqnarray*}
 \xi_0(x)&=&\xi  \\
\xi_1(x)&=&\sqrt x\, \psi - \xi /2   \\
\xi_2(x)&=&- \sqrt x \,\psi - \xi /2.
\end{eqnarray*}
 Both $\xi(x)$ and $\psi(x)$  have
 radius of convergence $x_c =4^{1/3}/3$.
\end{Lemma}
 The last statement is obtained using the Stirling formula. Note that one also has the following closed form expressions:
\begin{eqnarray*}
 \xi(x) &=& \displaystyle 2 \sqrt{\frac x 3} \  \sin
\left( 
\frac{\arcsin\left((3x)^{3/2} / 2\right)}{3}
\right),\\
\psi(x) &=& \displaystyle
\cos\left(  
\frac{\arcsin\left( (3x)^{3/2} / 2\right)}{3}
\right).
\end{eqnarray*}

\begin{figure}[ht]
\begin{center}
\epsfxsize=63mm \epsfysize=6cm
\epsf{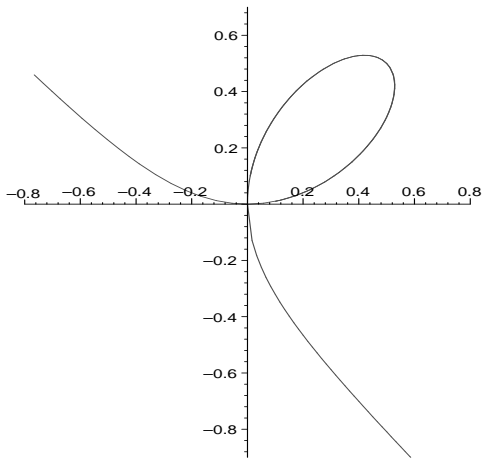}
\end{center}

\caption{The real values of  $\xi_i(x)$ for a real $x$.}
\label{descartes}
\end{figure}

We shall now consider the $\xi_i(x)$ as functions of a complex
variable $x$. 
We choose a determination of the  square root that coincides with
the usual determination on $\rs ^+ $:
$$ \sqrt{r e^{i \theta}} =
 \sqrt r e^{i \theta /2} \ \ \ \hbox{ for } \ \ \ -\pi  < \theta <
 \pi.$$
 Figure~\ref{descartes} shows the (real) values of the functions
 $\xi_i$ for real values of $x$. When $x$ is positive, the plots show,
 from bottom to top, $\xi_2(x)$, $\xi_0(x)$ and $\xi_1(x)$. When $x$
 is negative, the only real branch is $\xi_0$.

 The implicit function theorem implies that the
singularities of the functions $\xi_i$ are to be found among the
complex numbers $x$ such that the pair $(x,y)\equiv (x, \xi_i(x))$
satisfies $3y^2=x$ (and, of course, $x^3+y^3=xy$); that is,
\beq (x,y) \in \{ (0,0), (x_c,y_c), (j x_c, j^2y_c),  (j^2 x_c,
jy_c)\} \label{candidats}\eeq
where 
 $x_c=4^{1/3}/3$,  
$y_c=2^{1/3}/3$ and $j= \exp (2i\pi /3)$.
A more detailed investigation   gives the following result.
\begin{Lemma}
The singularities of the functions 
$\xi_i$ are given by:
$$Sing(\xi_0)= \{ x_c, j x_c , j^2 x_c\}, \ \ Sing(\xi_1)=\{0, x_c \}
\ \ \ \hbox{ and  }  \ \ \ Sing(\xi_2)=\{0,  jx_c,
j^2 x_c \}.$$
In particular, $\xi_2$ is not singular at $x_c$, as suggested by 
Figure~{\em\ref{descartes}}. All the above singularities are of the
square root type. 
\end{Lemma}
{\bf Proof.}
For each of the values of $x$ given by~\Ref{candidats}, we first have
to compute the values $\xi_i(x)$, $0 \le i \le 2$, in order to determine
which pairs $(x,\xi_i(x))$ are actually critical, and then, to check
the existence and nature of the singularity.

At $x=0$, all the $\xi_i(x)$ are zero. The explicit expansion of
Lemma~\ref{three-roots}, combined 
with the fact that $\psi$ and $\xi$ have a positive radius of
convergence, shows that only $\xi_1$ and $\xi_2$ are singular --- their
singularity being obviously of a square root type.

When $x=x_c$, factoring the polynomial $x^3+y^3-xy$ shows that the
multiset $\{\xi_0(x_c), \xi_1(x_c), \xi_2(x_c)\}$ equals
$\{y_c,y_c,-2y_c\}$. The functions $\xi_i(x)$ are real and continuous
on $[0,x_c]$, and the expansions of  Lemma~\ref{three-roots} show that
$\xi_0(x)$ and $\xi_1(x)$ are positive as $x\rightarrow 0^+$, while
$\xi_2(x)$ is negative. As the $\xi_i$ can only vanish at $0$, this
sign pattern must go on until $x_c$, so that
\beq \xi_0(x_c) = \xi_1(x_c)= y_c\ \ \ \  \hbox{ and } \ \ \  \xi_2(x_c)=
-2y_c.\label{values} \eeq
In view of~\Ref{candidats}, $x_c$ cannot be a singularity of
$\xi_2$. Now a local 
expansion of $x^3+y^3-xy$ around $(x_c,y_c)$ shows that as $x$
approaches $x_c^-$,
\beq \xi_{0,1}(x) =y_c \mp \frac 1 {\sqrt 3}
\sqrt{x_c-x}(1+o(1)),\label{square-root}\eeq
which confirms that the singularity of $\xi_0$ and $\xi_1$ at $x_c$ is
of the square root type.

 The values~\Ref{values} 
and Lemma~\ref{three-roots}
imply in particular that
 $$
\sqrt{x_c} \, \psi(x_c)=3y_c/2.
$$
Using this result, the fact that $\psi$ and $\xi$ are essentially
functions of $x^3$, and the values
$$\sqrt{j} = -j^2 \quad \hbox{ and } \quad  \sqrt{j^2} = -j,$$ 
we now compute
$$\begin{array}{lclcclcl}
\xi_0 (jx_c)  & = & j^2y_c  &&&\xi_0 (j^2 x_c)   &=  & jy_c   \\
\xi_1 (jx_c)   &=  &-2 j^2y_c  && &\xi_1 (j^2x_c)   &=  &  -2jy_c   \\
\xi_2 (jx_c)   & = & j^2y_c  && &\xi_2 (j^2x_c)   &=  &   jy_c  \\
\end{array}$$
so that $\xi_1$ cannot be singular at $jx_c$ or $j^2x_c$. Finally,
local expansions of  $x^3+y^3-xy$ confirm as above the existence of
square roots singularities of $\xi_0$ and $\xi_2$ at $jx_c$ and
$j^2x_c$.
\cqfd

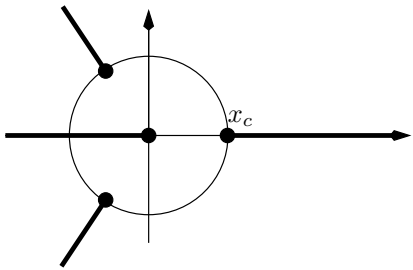
\begin{figure}[ht]
\begin{center}
\input{sing.pstex_t}
\end{center}
\caption{A domain on which all the functions $\xi_i$ are holomorphic.}
\label{plan}
\end{figure}

Each of the functions $\xi_i$ has a unique analytic continuation on any simply
connected domain avoiding $Sing(\xi_i)$, for instance the domain
obtained by removing  from $\cs$ the four 
half-lines of Figure \ref{plan}.

Observe that the series $\xi_1$ and $\xi_2$ can {\em also\/} be
 substituted for $y$ 
 in the functional equation~\Ref{equation-knight}; thus for $i \in \{0,1,2\}$, the following
 equation holds
\beq G(x)+G(\xi_i(x))=x^2\xi_i(x)^2, \label{func-eq} \eeq
at least  as an identity between power series
in $x$ or $\sqrt x$.

We end this subsection with a lemma that will be useful to build large
singularities from small ones.
\begin{Lemma}\label{grand}
Let $x \in \cs$, $x \not = 0$. Then one of the roots of $x^3+y^3-xy=0$
has modulus larger than $|x|$.
\end{Lemma}
{\bf Proof.} Let $y_0, y_1, y_2$ denote the three roots, and
assume none of then has modulus larger than $|x|$. The relation
$y_0y_1y_2=-x^3$ forces $|y_0|=|y_1|=|y_2|=|x|$. Then, the relation
$y_0+y_1+y_2=0$ implies $\{y_1, y_2 \}= \{jy_0 , j^2 y_0\}$. Finally,
the relation $y_0y_1+y_0y_2+y_1y_2=-x$ yields $x=0$. 
\qed

\subsection{The series $G(x)$ is not D-finite}
We now turn our attention to the series $G(x)$ defined by~\Ref{G-def}.
\begin{Proposition}\label{radius}
The series $G(x)$ has radius of convergence $x_c=4^{1/3}/3$. It is
singular at $x_c$, 
and, as $x \rightarrow x_c^-$, 
$$G(x) = A - B \sqrt{1-x/x_c} (1+o(1)),$$
where $A$ and $B$ are non-zero real numbers.
\end{Proposition}
{\bf Proof.} We start from  the functional equation defining $G$:
$$G(x)+G(\xi_0(x))=x^2 \xi_0(x)^2.$$
Recall that $G(x)$ and $\xi_0(x)$ have nonnegative coefficients
($G(x)$ counts walks, and $\xi_0(x)\equiv \xi(x)$ is given explicitly
by~\Ref{xi-def}). Hence  $G(\xi_0(x))$ and $x^2 \xi_0(x)^2$ also have 
nonnegative coefficients. 
The radius of convergence of $\xi_0(x)$ being
$x_c$, the series  $G(x)$ and  $G(\xi_0(x))$ have radius {\em at
least\/} $x_c$. 
The fact that $G(x)$ has radius at least $x_c$  can also be directly
derived from the upper bound~\Ref{bound}. 

For $|x|< x_c$, one has $|\xi_0(x)|< \xi_0(x_c)=y_c$, and $y_c$ is
smaller than $x_c$. This 
implies that the above functional equation also
holds as an identity between analytic functions of $x$ in the disk
$|x|< x_c$. As $x$ approaches $x_c$ inside this disk, a local
expansion gives, thanks to~\Ref{square-root},
$$G(x)= x_c^2y_c^2-G(y_c) -\frac 1 {\sqrt 3} \sqrt{x_c-x}\left(
2x_c^2y_c-G'(y_c)\right) + O(x_c-x).$$
The upper bound~\Ref{bound},
combined with~\Ref{G-def}, yields
$$G'(y_c) \le 3 \sum_{i \ge 1} (i+1) {{3i-2}\choose {i-1}} y_c^{3i+2}.$$
 This sum can be evaluated numerically, and is found to be smaller
than $0.16$. Comparing with  $2x_c^2y_c\approx 0.23$ gives the
announced result.
\cqfd

In view of Proposition~\ref{radius},
Eq.~\Ref{func-eq} now holds as an identity between
{\em functions\/} of $x$, as long as $|x|<x_c$ and $|\xi_i(x)|<x_c$ (and $x
\not \in \rs ^-$ if $i =1$ or $2$).

\begin{Proposition} \label{preuve-finale}
The series  $G(x)$ that counts knight walks ending on the $x$-axis is
not D-finite. 
\end{Proposition}
{\bf Proof.} 
Assume $G$ is D-finite. Then, it has a finite number of singularities,
and has a unique analytic continuation on any simply connected domain
avoiding these singularities.
 For $i \in \{0,1,2\}$, the series $G(\xi_i)$ are also D-finite (since
$\xi_i$ is algebraic).  Moreover, by analytic continuation of \Ref{func-eq},
\beq G(x)+G(\xi_i(x))=x^2\xi_i(x)^2 \label{cont} \eeq
for all $x$ in any domain  
where all our
functions are analytically defined. 
 Let us consider the
 identity~\Ref{cont} for $i=2$.
As $x$ approaches $x_c$, the function $G(x)$ becomes singular, while
 $\xi_2(x)$ does not: this shows that 
$\xi_2(x_c)=-2y_c$ is a singularity of $G$. 

Now let $x_s$ be a singularity of $G$ of maximal
modulus. 
According to Lemma~\ref{grand}, there exists $i$ such that 
$|\xi_i(x_s)|>|x_s|$. By assumption, $G(x)$ is singular at $x_s$.
But  $|x_s| \ge 2y_c >x_c$, so that $x_s$ is not a
singularity of $\xi_i$. 
Hence \Ref{cont} implies that $G$ is singular at $\xi_i(x_s)$. But
$|\xi_i(x_s)|>|x_s|$, which contradicts the maximality  of $|x_s|$.
\cqfd

\noindent{\bf Remarks} \\
1. The only property of D-finite series we have really used is the
fact that these series have  finitely many singularities. Hence,
 what we have actually proved is that {\em $G$ cannot have finitely many
singularities.} \\
2. The principle of the proof can be applied to other functional
equations of the same type. We apply it in the next section to a
(minor) variation of the knight walks.\\
3. As mentioned at the beginning of this section, the same walk
   problem was studied by Gessel in~\cite{gessel-proba}. This
   reference contains a formula analogous to~\Ref{eq:iterated}, and
   the statement that $G(x)$ has radius of convergence $4^{1/3}/3$.

\bigskip

There remains to prove Proposition~\ref{propo-diagonal}: the \gf \ for knight walks
ending on the diagonal is not D-finite. 

\noindent
{\bf Proof of Proposition~\ref{propo-diagonal}.} The length \gf \ for
walks ending 
on the diagonal is the coefficient of $x^0$ in the series
$t^{-2}Q(tx,t\bx)$, where $\bx =1/x$. This series is a power series in
$t$ whose coefficients are Laurent polynomials in $x$. Let us denote
$u=x^3$ and $\bu = 1/u$. Then, by Eq.~\Ref{equation-knight},
\beq t^2Q(tx,t\bx)=
\frac{t^4-S(t^3u)-S(t^3\bu)}{1-t(u+\bu)},\label{diagonale}\eeq 
where the series $S(z)$ is defined by
$S(z^3)=G(z).$
Let us convert $1/(1-t(u+\bu))$ in partial fractions of $u$. We obtain
$$\frac 1 {1-t(u+\bu)} = \frac 1 {\sqrt{1-4t^2}}\left(
\frac 1{1-u U(t)} +\frac 1{1-\bu U(t)} -1 \right),$$
where 
$$U(t)= \frac{1-\sqrt{1-4t^2}}{2t}.$$
We can now extract the constant term in $x$ 
(or in $u$)
from~\Ref{diagonale}:
\begin{eqnarray*}
[x^0] t^2Q(tx,t\bx)&=& \frac 1 {\sqrt{1-4t^2}}\left(
t^4 - [u^0] \frac {S(t^3u)}{1-\bu U(t)} 
 - [u^0] \frac {S(t^3\bu)}{1-u U(t)} \right)\\
&=&\frac 1 {\sqrt{1-4t^2}}\left(
t^4 -2S(t^3U(t))\right).
\end{eqnarray*}
Assume this series is D-finite: then the series $D(t)$ defined by
$D(t)=S(t ^3U(t))$ is D-finite in $t$ 
too. Let $T(s)$ be  the unique power series in $s$ such that $T(0)=0$,
 $T'(0)=1$
and $T^6-s^4T^2+s^8=0.$ Then the fact that D-finite series are stable
by any algebraic substitution tells us that $D(T(s))=S(T^3U(T))$ is
D-finite. But 
$T^3U(T)=s^4$, so that $S(s^4)$ itself, and hence $S(s)$ and $G(s)$,
are D-finite too, which we have proved to be false.
\cqfd

\section{A link with multidimensional linear recurrences with constant
coefficients}\label{section-recurrences}
In~\cite{bousquet-petkovsek}, we considered $d$-dimensional sequences
of complex 
numbers, denoted $a_{\bs{n}}=a_{n_1,n_2, \ldots n_d}$,
defined by recurrence relations of the following form:
\begin{equation}
\label{eqn}
a_{\bs{n}} =  \sum_{{\bs{h}}\in H} c_{\bs{h}} a_{{\bs{n}} + {\bs{h}}}\,\ \ \ \ \hbox{ for } {\bm{n}} \ge
{\bm{s}}, 
\end{equation}
\noindent 
where  $H = \{{\bm{h}}_1, {\bm{h}}_2, \ldots, {\bm{h}}_k\} 
\subseteq \zs^d$ is the set of {\em shifts\/},  $\left(c_{\bs{h}}\right)_{\bs{h} \in H}$ are given nonzero
constants, and ${\bm{s}} \in \ns^d$ is the 
{\em starting point\/} satisfying ${\bm{s}} + H \subseteq \ns^d$. 
We think of the ${\bm{h}}_i$ as having mostly (but not necessarily
only) negative coordinates, and of the
point ${\bm{n}}$ as depending on the points ${\bm{n}} + {\bm{h}}_1, 
{\bm{n}} + {\bm{h}}_2, \ldots, {\bm{n}} + {\bm{h}}_k$ as far as the value of
$a_{\bs{n}}$ is concerned. 
 A given function $\varphi$ specifies the initial conditions:
\begin{equation}
\label{init}
a_{\bs{n}} = \varphi({\bm{n}})\ \ \ \ \hbox{ for } {\bm{n}} \ge
{\bm{0}},\ 
{\bm{n}}\not\ge{\bm{s}}\,.
\end{equation}
   The convex hull of the set $H$ is assumed not to intersect the first
orthant  ({\em i.e.\/},  $\bm n \ge \bm 0$). This condition, as shown
in~\cite{bousquet-petkovsek}, guarantees that the numbers $a_{\bs{n}}$
can be computed recursively using~\Ref{eqn}. 

The enumeration of walks in a quadrant fits
exactly in this
framework, with $d=3$. Indeed, denoting $Q_{i,j}(n)$ the number of
walks that start
from a given point $(i_0,j_0)$, end at $(i,j)$ and
have length $n$, we have
$$Q_{i,j}(n) = \left\{
\begin{array}{ll}
1 & \hbox{if } (i,j,n)=(i_0,j_0,0) \\
0  & \hbox{if } i<0 \hbox{ or } j<0 \hbox{ or } n<0 \\
\displaystyle
\sum_{(h,k)\in\frak S} Q_{i-h,j-k}(n-1) & \hbox{otherwise,}
\end{array}
 \right.
$$
where $\frak S$ is the set of steps. A translation of the indices $i$
and $j$ transforms this recursion into one of the above type, with
$d=3$ and $H=\{(-h,-k,-1) \hbox{ such that } (h,k) \in \frak S\}$.

The paper \cite{bousquet-petkovsek} mostly dealt with the algebraic
nature of the generating function of the  
solution of such recurrences, and we started a classification, based on 
the {\em apex\/} of the recurrence,
defined as 
 the componentwise maximum of the points in
$H \cup \{\bm{0}\}$. We proved that when the initial conditions have
rational generating  
functions and the apex is $\bm{0}$, the generating function of the
solution is rational. Next,
when the initial conditions have algebraic generating 
functions and the apex has at most one positive coordinate, the 
generating function of the solution is algebraic.

When the apex has two positive coordinates, and $d=3$, our study of
walks in a quadrant shows that the solution 
might or might not be
D-finite. For 2-dimensional sequences, the ``simplest'' example with
apex $(1,1)$ was introduced in~\cite{petkovsek-toronto}: 
For $i,j \ge 0$,  let
\beq a_{i,j}= \left\{ 
\begin{array}{ll}
a_{i+1,j-2}+a_{i-2,j+1} & \hbox{ if } i,j \ge 2,\\
1 & \hbox{ otherwise. }
\end{array} \right. \label{rec2} \eeq
This
recurrence is obviously closely connected to the knight walk. Again, it has a
unique solution whose terms can be computed 
inductively. The first few values
are given in the following array.
$$\begin{array} {c|cccccccc}
j &&&&&&& \\
\uparrow &&&&&&& \\
6 & 1	&	1	&	5	&	7	& \cdot &
\cdot & \cdot & \\

5 & 1 & 1 & 3 & 5 & 10 & 14 & \cdot & \\

4 & 1 & 1 & 3 & 4 &  6 & 10 & \cdot & \\

3 & 1 & 1 & 2 & 2 &  4 &  5 &  7 & \\

2 & 1 & 1 & 2 & 2 &  3 &  3 &  5 & \\

1 & 1 & 1 & 1 & 1 &  1 &  1 &  1 & \\

0 & 1 & 1 & 1 & 1 &  1 &  1 &  1  \\ \hline
 & 0 & 1 & 2 & 3 & 4 & 5 & 6 & \rightarrow \ i \\
\end{array}$$

\noindent Defining the generating function
$$A(x,y)=\sum_{i,j \ge 2} a_{i,j} x^{i-2} y^{j-2},$$ we obtain, by summing the
recurrence relation over $i, j \ge 2$, the following functional
equation:
\beq ( xy-x^3-y^3)A(x,y)=R(x,y)-F(x)-F(y)\label{cavalier} \eeq
where 
$$R(x,y) = xy\left(\frac{1+y}{1-x}+\frac{1+x}{1-y}\right) \ \ \ \ \hbox{ and
}\ \ \ \ 
 F(x)=\sum_{i\ge 2} a_{i,2} x^{i+1}=x^3A(x,0).$$
We have used the symmetry of the problem in $i$ and $j$.
It was first proved in~\cite{petkovsek-toronto} that  $F(x)$
and $A(x,y)$ are irrational. Then, we claimed
in~\cite[p.~74]{bousquet-petkovsek} that they 
are even not D-finite, but
without giving a proof. The tools developed above for the knight walk
apply perfectly to this problem.
\begin{Proposition}
%
The \gf \ $A(x,y)$ of the bivariate sequence $a_{i,j}$ defined by the
recurrence relation~{\em(\ref{rec2})} is not D-finite. Nor is the series
$A(x,0)$.
\end{Proposition}
{\bf Proof.}
The argument is very close to the knight one. The  kernel
method first gives 
\beq  F(x)+F(\xi_i(x))
=R(x,\xi_i(x)),\label{main2} \eeq
for any of the three roots $\xi_i$ of the kernel.

The only difference with the knight treatment comes from the fact
that  $R(x, \xi_i(x))$ might have --- and indeed, has --- more
singularities than $\xi_i$. They can be determined exactly, but
we shall only use the following obvious information
\beq Sing(R(x, \xi_i(x))) \subset Sing(\xi_i) \cup \{1\} \cup \{x :
x^3-x+1=0\}.\label{sing-K} \eeq
In particular, all singularities of $R(x, \xi_i(x))$ have modulus at most
 $1.33$ (an upper bound for the modulus of the largest root of $x^3-x+1$). 
We start from Eq.~\Ref{main2}, in the case $i=0$. 
As  $R(x,\xi_0(x))$ has radius $x_c$, we
find again that the radius of $F(x)$ is at least $x_c$. Moreover,
comparing the recurrence relations~\Ref{rec1} and~\Ref{rec2}, and the
corresponding tables, shows that for $i,j \ge 2$, one has $a_{i,j}\ge
 Q_{i-2,j-2}$. This implies that the radius of $F(x)$ is bounded from
above by the radius of $G(x)$, which was proved to be exactly
$x_c$. Hence $F(x)$ has radius $x_c$.

Assume $F(x)$ is D-finite. As in the proof of
Proposition~\ref{preuve-finale}, we first construct large
singularities of $F$ (here, large means larger than $1.33$).
Let us start from $x_0 = x_c=4^{1/3}/3\approx 0.53$,
which {\em is\/} a
singularity of $F$
(by Pringsheim theorem, a series with nonnegative coefficients is
singular at its radius of convergence). 
As it is not a singularity of $R(x, \xi_2(x))$,
Eq.~\Ref{main2} implies that
$x_1:=\xi_2(x_c)=-2y_c\approx-0.84$ is singular for $F$.

But  $x_1$ is singular for none of the $R(x,
\xi_i(x))$. Moreover,
$$\{ \xi_0(x_1), \xi_1(x_1), \xi_2(x_1) \}= \{ x_0, -0.26...\pm
1.02... i\}.$$
Using the same trick as above, we see that $x_2\approx-0.26-1.02i$ is
singular for $F$.

One more step: Among the $\xi_i(x_2)$, one is
$x_3\approx 0.92-1.02 i$, 
which has modulus larger than $1.33$. 
As $x_2$ is not singular for any of the $R(x, \xi_i(x))$, the $\xi_i(x_2)$
are singularities of $F$. 
In particular, $x_3$ is singular for $F$.

We conclude as above, by considering the largest singularity $x_s$ of $F$
(in modulus), showing that one of the $\xi_i(x_s)$ is singular for $F$
{\em and\/} larger than $|x_s|$ in modulus.

\cqfd

\noindent {\bf Acknowledgements.} We are grateful to both referees for
pointing out interesting references and simplifications in the proof
of Lemma~\ref{three-roots}.

\end{document}

%% file: bijection.pstex_t
\begin{picture}(0,0)%
\includegraphics{bijection.pstex}%
\end{picture}%
\setlength{\unitlength}{4144sp}%
\begingroup\makeatletter\ifx\SetFigFont\undefined%
\gdef\SetFigFont#1#2#3#4#5{%
  \reset@font\fontsize{#1}{#2pt}%
  \fontfamily{#3}\fontseries{#4}\fontshape{#5}%
  \selectfont}%
\fi\endgroup%
\begin{picture}(6087,2634)(451,-2008)
\put(991,-1591){\makebox(0,0)[lb]{\smash{\SetFigFont{10}{12.0}{\familydefault}{\mddefault}{\updefault}\special{ps: gsave 0 0 0 setrgbcolor}$i_0$\special{ps: grestore}}}}
\put(5131,-1816){\makebox(0,0)[lb]{\smash{\SetFigFont{10}{12.0}{\familydefault}{\mddefault}{\updefault}\special{ps: gsave 0 0 0 setrgbcolor}$\bar s_2$\special{ps: grestore}}}}
\put(451,-61){\makebox(0,0)[lb]{\smash{\SetFigFont{10}{12.0}{\familydefault}{\mddefault}{\updefault}\special{ps: gsave 0 0 0 setrgbcolor}$j_0$\special{ps: grestore}}}}
\put(3556,-106){\makebox(0,0)[lb]{\smash{\SetFigFont{10}{12.0}{\familydefault}{\mddefault}{\updefault}\special{ps: gsave 0 0 0 setrgbcolor}$j_0$\special{ps: grestore}}}}
\put(496,-556){\makebox(0,0)[lb]{\smash{\SetFigFont{10}{12.0}{\familydefault}{\mddefault}{\updefault}\special{ps: gsave 0 0 0 setrgbcolor}$4$\special{ps: grestore}}}}
\put(541,-1006){\makebox(0,0)[lb]{\smash{\SetFigFont{10}{12.0}{\familydefault}{\mddefault}{\updefault}\special{ps: gsave 0 0 0 setrgbcolor}$2$\special{ps: grestore}}}}
\put(541,-1231){\makebox(0,0)[lb]{\smash{\SetFigFont{10}{12.0}{\familydefault}{\mddefault}{\updefault}\special{ps: gsave 0 0 0 setrgbcolor}$1$\special{ps: grestore}}}}
\put(2296,-1096){\makebox(0,0)[lb]{\smash{\SetFigFont{10}{12.0}{\familydefault}{\mddefault}{\updefault}\special{ps: gsave 0 0 0 setrgbcolor}$s_2$\special{ps: grestore}}}}
\put(1711,-1366){\makebox(0,0)[lb]{\smash{\SetFigFont{10}{12.0}{\familydefault}{\mddefault}{\updefault}\special{ps: gsave 0 0 0 setrgbcolor}$s_1$\special{ps: grestore}}}}
\put(4591,-1591){\makebox(0,0)[lb]{\smash{\SetFigFont{10}{12.0}{\familydefault}{\mddefault}{\updefault}\special{ps: gsave 0 0 0 setrgbcolor}$\bar s_1$\special{ps: grestore}}}}
\put(4141,-1546){\makebox(0,0)[lb]{\smash{\SetFigFont{10}{12.0}{\familydefault}{\mddefault}{\updefault}\special{ps: gsave 0 0 0 setrgbcolor}$i_0$\special{ps: grestore}}}}
\end{picture}

%% file: sing.pstex_t
\begin{picture}(0,0)%
\epsfig{file=sing.pstex}%
\end{picture}%
\setlength{\unitlength}{3947sp}%
\begingroup\makeatletter\ifx\SetFigFont\undefined%
\gdef\SetFigFont#1#2#3#4#5{%
  \reset@font\fontsize{#1}{#2pt}%
  \fontfamily{#3}\fontseries{#4}\fontshape{#5}%
  \selectfont}%
\fi\endgroup%
\begin{picture}(2610,1698)(1318,-1141)
\put(2746,-196){\makebox(0,0)[lb]{\smash{\SetFigFont{10}{12.0}{\rmdefault}{\mddefault}{\updefault}$x_c$}}}
\end{picture}